\catcode`\à=\active \def à{\` a}
\catcode`\â=\active \def â{\^ a}
\catcode`\é=\active \def é{\' e}
\catcode`\è=\active \def è{\` e}
\catcode`\ê=\active \def ê{\^ e}
\catcode`\ë=\active \def ë{\" e}
\catcode`\ù=\active \def ù{\` u}
\catcode`\û=\active \def û{\^ u}
\catcode`\ô=\active \def ô{\^ o}
\catcode`\ç=\active \def ç{\c c}
\catcode`\î=\active \def î{\^ \i}
\catcode`\ï=\active \def ï{\" \i}

\magnification=1100

\vsize=215mm

\baselineskip 13pt

\font\zazouete = cmb10 at 9pt

\font\sc=cmcsc10 at 10pt

\font\sccp=cmcsc10 at 9pt

\font\pg=cmss10 scaled 950

\font\tfp=cmssbx10 scaled 1200

\font\tfpp=cmssbx10 scaled 1000

%\NoBlackBoxes

%\nopagenumbers

\font\mathb=msbm10 scaled 1000

\let\Bbb=\mathbbb

\input amssym

\input epsf

\font\ib = cmbxsl10

\centerline{\zazouete COHOMOLOGIE DE DOLBEAULT FEUILLET\'EE}
%\smallskip
\centerline{\zazouete DU FEUILLETAGE COMPLEXE AFFINE DE REEB}

\smallskip

\centerline{par}

\smallskip

\centerline{Rochdi {\sc Ben Charrada} \& Aziz {\sc El Kacimi Alaoui}}

\footnote{}{{\it Mathematics Subject Classification } : 32W05,
32G05, 32Q58, 58A30

}

\medskip
\centerline{(Septembre  2019)}

\vskip0.3cm
\noindent {\parindent=0.7cm\narrower {\pg {\tfpp R\'esum\'e.}
Soit ${\cal F}$ le  feuilletage complexe affine de Reeb de dimension $1$ sur  la vari\'et\'e
de Hopf \hskip0.1cm ${\Bbb S}^{n+1}\times {\Bbb S}^1$. On montre que sa  cohomologie   de Dolbeault feuillet\'ee
en degré $1$  est
isomorphe \`a ${\Bbb C}$ en en exhibant explicitement un g\'en\'erateur. On voit apparaître  ainsi toutes les obstructions à résoudre
le $\overline{\partial }_{\cal F}$ le long des feuilles sur $({\Bbb S}^{n+1}\times {\Bbb S}^1,{\cal F})$.
\par }}

\vskip0.5cm

\medskip
{\tfp 1. Premi\`eres d\'efinitions}
\medskip

\noindent Soit $M$ une vari\'et\'e diff\'erentiable (de classe $C^\infty $) de dimension $m+n$. On suppose, pour simplifier, qu'elle est connexe
et qu'elle poss\`ede toutes les bonnes propri\'et\'es dont on aurait \'eventuellement besoin (paracompacit\'e...).

\newif\ifpremierepage\premierepagefalse
\def\makeheadline{%
\ifpremierepage\global\premierepagefalse \else \vbox to
0pt{\vskip-22.5pt \line{\vbox to 8.5pt{}\ifodd\pageno \textedroite
\else \textegauche \fi } \vskip4pt
%\hrule\vss
}\nointerlineskip \fi}
\def\boxit#1#2{\setbox1=\hbox{\kern#1{#2}\kern#1}%
\dimen1=\ht1 \advance\dimen1 by #1 \dimen2=\dp1 \advance\dimen2 by
#1
\setbox1=\hbox{\vrule height\dimen1 depth\dimen2\box1\vrule width 1mm}%
\setbox1=\vbox{\hrule\box1\hrule height 0.5mm depth 0.5mm}%
\advance\dimen1 by .4pt \ht1=\dimen1 \advance\dimen2 by .4pt
\dp1=\dimen2 \box1\relax}

\gdef\textegauche{{ } \hfil {\sccp R. Ben Charrada \& A. El Kacimi Alaoui }\hfil }
\gdef\textedroite{{ }\hfil {\pg Cohomologie de Dolbeaut feuillet\'ee du feuilletage affine complexe de Reeb}\hfil }
\premierepagetrue

\smallskip
\noindent {\tfpp 1.1. D\'efinion.} {\it  Un
{\ib feuilletage} ${\cal F}$ de {\ib codimension} $n$ sur $M$ est
donn\'e par un recouvrement ouvert ${\cal U}=\{ U_i\}_{i\in I}$ et,
pour tout $i$, d'un diff\'eomorphisme ${\Bbb
R}^m\times {\Bbb R}^n\buildrel {\varphi_i}\over \longrightarrow U_i$ tel que, sur toute intersection non vide
  $U_i\cap U_j$, le diff\'eomorphisme de changement de coordonn\'ees :
  $$\varphi_j^{-1}\circ \varphi_i: (z,t)\in \varphi_i^{-1}(U_i\cap U_j)
\longrightarrow(z',t')\in \varphi_j^{-1}(U_i\cap U_j)$$
soit de la forme $z'=\varphi_{ij}(z,t)$ et $t'=\gamma_{ij}(t)$.}
\smallskip
La vari\'et\'e  $M$ est ainsi d\'ecompos\'ee en sous-vari\'et\'es connexes de dimension $m$.
Chacune d'elles est appel\'ee {\it feuille} de ${\cal F}$. On note $\tau $ le fibr\'e tangent \`a ${\cal F}$ ;
il est constitu\'e de tous les vecteurs  tangents aux feuilles.
Les sections de $\tau $   forment un module
$\frak{X}({\cal F})$ sur l'anneau $C^\infty (M)$ des fonctions  sur $M$. Si $X$ et $Y$ sont deux \'el\'ements de
$\frak{X}({\cal F})$, par le th\'eor\`eme de Frobenius [Ca], le crochet $[X,Y]$ est encore un \'el\'ement de $\frak{X}({\cal F})$.
\smallskip
Le quotient $\nu {\cal F}=TM/\tau $ est le {\it fibr\'e normal} \`a ${\cal F}$ ; on peut le r\'ealiser dans $TM$ par
le choix d'un suppl\'ementaire $\nu $ de $\tau $. On a ainsi une d\'ecomposition en somme directe $TM=\tau \oplus \nu $.
Celle-ci donne une d\'ecomposition du complexifi\'e du fibr\'e $\Lambda^\ell T^\ast M\otimes {\Bbb C}$ des $\ell $-formes ext\'erieures :
$$\Lambda^\ell T^\ast M\otimes {\Bbb C}=\bigoplus_{s+r=\ell }\Lambda^{sr}\leqno{(1)}$$
o\`u $\Lambda^{sr}$ est le fibr\'e dont les sections globales sont les $\ell $-formes complexes $\alpha $ de type $(s,r)$ {\it i.e.} celles qui s'\'ecrivent localement :
$$\alpha =\sum_{1\leq i_1<\cdots <i_s\leq m\atop 1\leq j_1<\cdots <j_r\leq }f_{i_1,\cdots ,i_s, j_1,\cdots ,j_r}(z,t)dt_{i_1}\wedge \cdots \wedge dt_{i_s}
\wedge dz_{j_1}\wedge \cdots \wedge dz_{j_r}\leqno{(2)}$$
o\`u les $f_{i_1,\cdots ,i_s, j_1,\cdots ,j_r}$ sont des fonctions continues et $C^\infty $ en $z=(z_1,\cdots ,z_m)$.
L'ensemble $A^{sr}(M)$ de ces formes diff\'erentielles est un module sur l'anneau $A(M)=C^{0,\infty }(M,{\Bbb C})$ des fonctions complexes
sur $M$ ($C^0$ en $(z,t)$ mais $C^\infty $ en $z$).
\smallskip

Dans toute la suite, on se restreindra au cas $s=0$. On note alors $A_{\cal F}^r(M)$ l'espace $A^{0r}(M)$ et on considère l'opérateur
$d_{\cal F}$ qui à la forme $\alpha \in A_{\cal F}^r(M)$ associe la forme
$d_{\cal F}\alpha \in A_{\cal F}^{r+1}(M)$ dont l'évaluation $d_{\cal F}\alpha (X_1,\cdots ,X_{r+1})$ sur les $r+1$ champs de vecteurs $X_1,\cdots ,X_{r+1}$ tangents à
${\cal F}$ est donnée  par :
$$\sum_{i=1}^{r+1}(-1)^iX_i\cdot \alpha (X_1,\cdots ,\widehat X_i,\cdots ,X_{r+1})+
\sum_{i<j}(-1)^{i+j}\alpha ([X_i,X_j],\cdots ,\widehat X_i,\cdots \widehat X_j,\cdots X_{r+1}).\leqno{(3)}$$

Cet opérateur est de carré nul et donne alors un complexe différentiel :
$$0\longrightarrow A_{\cal F}^0(M)\buildrel {d_{\cal F}} \over \longrightarrow A_{\cal F}^1(M)\buildrel {d_{\cal F}} \over \longrightarrow \cdots \buildrel {d_{\cal F}} \over \longrightarrow A_{\cal F}^{m-1}(M)
\buildrel {d_{\cal F}} \over \longrightarrow A_{\cal F}^m(M)\longrightarrow 0\leqno{(4)}$$
appel\'e {\it complexe de de Rham feuillet\'e} de ${\cal F}$.
Son homologie en degr\'e $r$ sera not\'ee $H_{\cal F}^r(M)$ et appel\'ee {\it cohomologie feuillet\'ee} de ${\cal F}$.
Elle co\"\i ncide avec la cohomologie de de Rham de $M$ lorsque la dimension des feuilles est celle de la vari\'et\'e, c'est-\`a-dire lorsque
il n'y a qu'une seule feuille, la vari\'et\'e $M$ elle-m\^eme.

\smallskip
On se donne maintenant une vari\'et\'e  $M$ comme avant, qu'on suppose   de dimension $2m+n$
et munie d'un feuilletage ${\cal F}$ de codimension $n$ (et donc de dimension $2m$).
\smallskip
\noindent {\tfpp 1.2. D\'efinition.} {\it  On dira que ${\cal F}$
est {\ib complexe}  s'il existe un recouvrement ouvert $\{
U_i\} $ de $M$ et des diff\'eomorphismes $\phi_i : \Omega_i\times
{\cal O}_i\longrightarrow U_i$, o\`u $\Omega_i$ est un ouvert de
$\Bbb C^m$ et ${\cal O}_i$ un ouvert de $\Bbb   R^n$ tels que les
changements de coordonn\'ees :
$$\phi_{ij}=\phi_j^{-1}\circ \phi_i:\phi_i^{-1}(U_i\cap U_j)\longrightarrow
\phi_j^{-1}(U_i\cap U_j)$$
soient de la forme $(z',t')=\left( \phi_{ij}^1(z,t),\phi_{ij}^2(t)\right) $
avec $\phi_{ij}^1(z,t)$
holomorphe en $z$ pour $t$ fix\'e.}
\smallskip
Chaque feuille de ${\cal F}$ est une vari\'et\'e analytique complexe de dimension $m$.
La notion
de feuilletage complexe g\'en\'eralise celle
de feuilletage holomorphe sur une vari\'et\'e analytique
complexe.
\smallskip
La donn\'ee d'un feuilletage complexe ${\cal F}$ sur une
vari\'et\'e $M$ sera repr\'esent\'ee par le  couple $(M,{\cal F})$.

Soient $(M,{\cal F})$ et $(M',{\cal F}')$ deux feuilletages
complexes. On appelle {\it morphisme} de $(M,{\cal F})$ vers
$(M',{\cal F}')$ toute application $f:M\longrightarrow M'$, de
classe $C^\infty $ et telle que l'image de toute feuille $F$ de
${\cal F}$ est contenue dans une feuille $F'$ de ${\cal F}'$ et
 l'application $f:F\longrightarrow F'$ est holomorphe.

Un morphisme $f:(M,{\cal F})\longrightarrow (M',{\cal
F}')$ est un {\it isomorphisme de feuilletages complexes} si c'est
un diff\'eo\-morphisme qui est un biholomorphisme sur les feuilles.
On dira que deux feuilletages complexes ${\cal F}$ et ${\cal F}'$
sur $M$ sont {\it conjugu\'es} (ou dans la m\^eme {\it classe de
conjugaison}) s'il existe un isomorphisme $f:(M,{\cal
F})\longrightarrow (M,{\cal F}')$.
\smallskip
L'ensemble des automorphismes de ${\cal F}$ est un groupe qu'on
notera $G({\cal F})$. On peut remarquer qu'une feuille de ${\cal F}$
qui n'est biholomorphiquement \'equivalente \`a aucune autre feuille
de ${\cal F}$   est fix\'ee par tout le groupe $G({\cal F})$.
\smallskip
\noindent {\tfpp 1.3. Exemples}

\smallskip
\item{i)} Toute vari\'et\'e analytique complexe de dimension $m$ est un
feuilletage complexe de dimension $m$. Le groupe des
automorphismes de ce feuilletage est r\'eduit \`a celui des
automorphismes de la vari\'et\'e complexe.
\smallskip
\item{ii)} Tout feuilletage holomorphe au sens usuel sur une vari\'et\'e complexe
est un feuilletage complexe sur la vari\'et\'e r\'eelle
sous-jacente.
\smallskip
\item{iii)} Soit $M$ un ouvert de $\Bbb C^m\times \Bbb R^n$. Pour tout $t\in \Bbb R^n$, on note $M^t$
l'ensemble :
$$\{ z\in \Bbb C^m:(z,t)\in M \} .$$
$M^t$ est un ouvert
de $\Bbb C^m$   appel\'e  {\it section} de $M$ suivant $t$. Les
sections   sont les feuilles d'un feuilletage complexe ${\cal F}$ de
dimension $m$   qu'on appellera {\it feuilletage complexe
canonique} de $M$.
\smallskip
\item{iv)} Soit $F$ une vari\'et\'e analytique complexe de dimension $m$. Toute fibration
localement triviale  $F\hookrightarrow M\longrightarrow B$ dont le
cocycle est \`a valeurs dans le groupe Aut$(F)$ des biholomorphismes
de la fibre $F$ est un feuilletage complexe ${\cal F}$ de
dimension $m$. Si la fibration est triviale
{\it i.e.} $M=F\times B$, on dira que $F$ est un {\it feuilletage
produit} : toutes les feuilles $F\times \{ t\} $ ont la m\^eme
structure complexe.
\smallskip
\item{v)} Supposons que ${\cal F}$ est un feuilletage
complexe sur $M=F\times B$ dont les feuilles sont les facteurs
$F\times \{ t\} $ mais que la structure complexe n'est pas
forc\'ement la m\^eme  sur toutes  les feuilles ; on dira alors que
${\cal F}$ est un {\it produit diff\'erentiable}.
\smallskip
\item{vi)} Soit ${\cal F}$ un feuilletage
orientable par surfaces  sur une vari\'et\'e $M$. On consid\`ere une
m\'etrique riemannienne $g$ sur le fibr\'e $T{\cal F}$. A tout
vecteur
 $u\in T_y{\cal F}$ on associe l'unique vecteur
$v\in T_y{\cal F}$ de m\^eme longueur que $u$ et tel que le rep\`ere
$(u,v)$ soit direct. En posant $Ju=v$ on d\'efinit ainsi une
structure presque complexe sur chaque feuille. La version \`a
param\`etre du th\'eor\`eme  d'int\'egrabilit\'e montre que cette
structure est en fait une {\it structure complexe} transversalement
(localement) diff\'erentiable sur ${\cal F}$. Ainsi tout feuilletage
orientable par surfaces
%admet une structure de
est un feuilletage
complexe de dimension $1$.
\medskip
{\tfp 2. La $\overline \partial _{\cal F}$-cohomologie}
\medskip
\noindent {\tfpp 2.1.} Soit $(M,{\cal F})$ un feuilletage complexe de dimension
$m$. On s'int\'eresse aux formes
feuillet\'ees qui,  dans un syst\`eme de coordonn\'ees
locales $(z,t)=(z_1,\ldots ,z_m,t_1,\ldots ,t_n)$, s'écrivent :
$$\alpha = \sum f_{j_1\ldots j_pk_1\ldots k_q}(z,t)dz_{j_1}\wedge
\ldots \wedge dz_{j_p} \wedge d\overline z_{k_1}\wedge \ldots \wedge
d\overline z_{k_q}$$ o\`u $f_{j_1\ldots j_pk_1\ldots k_q}$ est
une fonction continue   en $(z,t)$ mais $C^\infty $ en $z$. On les appelle {\it formes feuillet\'ees de type} $(p,q)$.
Elles forment un espace vectoriel qu'on notera $A_{\cal F}^{pq}(M)$ et qui est aussi un module sur
$A(M)$ (anneau des fonctions continues et $C^\infty $ en $z$). Ainsi, toute forme $\alpha \in A_{\cal F}^r(M)$
se d\'ecompose en une somme  $\alpha =\displaystyle \sum_{p+q=r}\alpha_{pq}$
o\`u $\alpha_{pq}$ est une forme feuillet\'ee de type $(p,q)$.
Ce qui donne la d\'ecomposition en somme directe :
$$A_{\cal F}^r(M)=\bigoplus_{p+q=r}A_{\cal F}^{pq}(M).\leqno{(5)}$$
On fixe l'entier $p\in \{ 0,1,\cdots ,m\} $. Alors l'{\it op\'erateur de Cauchy-Riemann} le long des feuilles   $\overline
\partial _{\cal F}:A^{pq}({\cal F})\longrightarrow A^{p,q+1}({\cal F})$
 s'\'ecrit localement :
$$\overline \partial _{\cal F}(\alpha (z,t)dz_{j_1}\wedge
\cdots \wedge dz_{j_p} \wedge d\overline z_{k_1}\wedge \ldots \wedge
d\overline z_{k_q})=$$
$$\sum_{s=1}^m{{\partial \alpha }\over {\partial \overline z_s}}(z,t)d\overline z_s\wedge
dz_{j_1}\wedge \ldots \wedge dz_{j_p} \wedge d\overline
z_{k_1}\wedge \ldots \wedge d\overline z_{k_q}$$ o\`u ${{\partial
\;}\over {\partial \overline z_s}}={{1}\over {2}}\big\{ {{\partial
\;}\over {\partial x_s}}+ i{{\partial \;}\over {\partial y_s}}\big\}
$ avec $z_s=x_s+iy_s$. On peut   v\'erifier facilement que cet op\'erateur est de carr\'e nul et qu'on  a un complexe diff\'erentiel :
$$0\longrightarrow A_{\cal F}^{p0}(M)\buildrel {\overline \partial _{\cal F}} \over \longrightarrow
A_{\cal F}^{p1}(M)\buildrel {\overline \partial _{\cal F}} \over
\longrightarrow \cdots \buildrel {\overline \partial _{\cal F}}
\over \longrightarrow A_{\cal F}^{p,m-1}(M) \buildrel {\overline \partial _{\cal F}}
\over \longrightarrow A_{\cal F}^{pm}(M) \longrightarrow  0$$
appellé {\it complexe du $\overline \partial_{\cal F}$} ou {\it complexe de Dolbeault feuillet\'e} de
$(M,{\cal F})$ ; son homologie, not\'ee $H_{\cal F}^{pq}(M)$, sera appel\'ee
la {\it $\overline \partial_{\cal F}$-cohomologie} ou {\it cohomologie de Dolbeault feuillet\'ee} de ${\cal F}$. Lorsque $M$ est une  vari\'et\'e
complexe munie du feuilletage dont la seule feuille est elle-m\^eme, alors
$H_{\cal F}^{pq}(M)$ n'est rien d'autre que sa cohomologie de Dolbeault usuelle.
%Le lecteur d\'esireux d'en savoir plus peut consulter [Ek1], [ES], [BC] ou [GT].
\medskip
\noindent {\tfpp 2.2.} Une $p$-forme feuillet\'ee $\omega $ est dite ${\cal F}$-{\it holomorphe} si elle est de type $(p,0)$
et v\'erifie $\overline \partial_{\cal F}\omega =0$. Le faisceau des germes de telles formes sera not\'e ${\cal O}_{\cal F}^p$ ; il admet
une r\'esolution fine :
$$0\longrightarrow {\cal O}_{\cal F}^p \hookrightarrow  {\cal A}_{\cal F}^{p0} \buildrel {\overline \partial _{\cal F}} \over \longrightarrow
{\cal A}_{\cal F}^{p1} \buildrel {\overline \partial _{\cal F}} \over \longrightarrow
\cdots \buildrel {\overline \partial _{\cal F}} \over \longrightarrow
{\cal A}_{\cal F}^{p,m-1} \buildrel {\overline \partial _{\cal F}} \over \longrightarrow
{\cal A}_{\cal F}^{pm}  \longrightarrow 0$$
o\`u ${\cal A}_{\cal F}^{pq} $ est le faisceau des germes des formes feuillet\'ees de type $(p,q)$. Comme
$A_{\cal F}^{pq}(M)$ est l'espace des sections globales du faisceau ${\cal A}_{\cal F}^{pq}$, on a un isomorphisme canonique :
$$H_{\cal F}^{pq}(M)\simeq H_{\cal F}^{q}(M,{\cal O}_{\cal F}^p).\leqno{(6)}$$
Les deux d\'efinitions permettent de faire des calculs suivant la nature des exemples et la mani\`ere dont ils sont d\'ecrits.
Nous verrons dans la suite comment cela se passe.
L'espace ${\cal H}_{\cal F}^p(M)$ des $p$-formes ${\cal F}$-holomorphes (sections globales du faisceau ${\cal O}_{\cal F}^p$)
sur $M$ est $H_{\cal F}^{p0}(M)$.

\smallskip

Divers  calculs  de la cohomologie feuillet\'ee de Dolbeault, et quelques-unes de leurs applications, sont donnés
dans [BC], [Ek1], [Ek2],  [ES], [GT], [S$\ell $].
\medskip

{\tfp 3. Fonctions ${\cal F}$-holomorphes en dimension 1}
\smallskip
\noindent Dans toute cette section, $M$ sera une vari\'et\'e
r\'eelle de dimension $2+n$ munie d'un feuilletage complexe ${\cal
F}$ de dimension $1$. L'espace
$A^{00}({\cal F})$ n'est rien d'autre que l'anneau $C^{0,\infty }(M)$ des fonctions complexes sur $M$ continues et $C^\infty $ le long des
feuilles qu'on a d\'ej\`a not\'e $A(M)$.
\smallskip
\noindent {\tfpp 3.1. La topologie $C^{0,\infty }$ sur
$A_{\cal F}^{0q}(M)$}
\smallskip

Soient $V$ un ouvert de ${\Bbb C}$ et $O$ un ouvert de ${\Bbb R}^n$ ;
les coordonn\'ees sur $V\times O$ seront not\'ees
$(z,t)=(z,t_1,\cdots ,t_n)$. Pour tout multi-indice $k=(k_1,k_2)$ ($k_1$ et $k_2$ sont des entiers naturels), on posera $\vert k\vert
=k_1+k_2$   et :
$$D^k={{\partial^{\vert k\vert }}\over {\partial z^{k_1}\partial
\overline z^{k_2}}}.$$
Pour tout $s\in \Bbb N$, tout compact $K$ de $\Omega \times O$  et
toute fonction $f:\Omega \times O\longrightarrow \Bbb
C$ on pose :
$$ N_{s,K}(f)=\max_{\vert k\vert \leq s}\left\{ \sup_K\left\vert D^k (f)\right\vert \right\}  .$$
Soit  $U$ un ouvert de la variété $M$
distingu\'e pour ${\cal F}$ et \'equivalent \`a $V\times O$ via un
isomorphisme $\varphi : V\times O\longrightarrow U$. Pour tout
compact $K$ contenu dans $U$ et toute fonction $f\in C^{0,\infty }(U)$
on pose $ N_{s,K}(f)=N_{s,K}(f\circ \varphi ) .$
\smallskip
Soient ${\cal U}=\{ (U,\varphi )\} $ un atlas d\'enombrable
d\'efinissant ${\cal F}$ et ${\cal C}=\{ C_n\} $ une suite de compacts, chacun contenu dans une carte de ${\cal
U}$,  recouvrant $M$ et telle que tout compact $K\subset M$ soit
recouvert par un nombre fini d'\'el\'ements de ${\cal C}$.
Consid\'erons une suite  croissante de compacts $K_n$ dont la
r\'eunion est \'egale \`a $M$. Pour tout $n\in \Bbb N^\ast $,
l'ensemble ${\cal C}_n$  des compacts de la famille ${\cal C}$ qui
intersectent  $K_n$ est fini. Pour tout $s\in \Bbb N$ et toute
fonction $f\in A(M)$, posons :
$$\vert \vert f\vert \vert_s^n=\sum_{C\in {\cal C}_n}N_{s,C}(f).$$
  La famille des semi-normes $\vert \vert \;\; \vert \vert_s^n$ (index\'ee
par $s\in \Bbb N$ et $n\in \Bbb N^\ast $) est filtrante et
s\'eparante ; elle permet de d\'efinir une distance sur $A(M)$
invariante par translations
$$\delta (f,g)=\sum_{s,n}{{1}\over {2^{s+n}}}\inf \left( 1, \vert \vert
f-g\vert \vert_s^n\right) .$$
Cette distance d\'efinit une topologie faisant de $C^{0,\infty }(M)$ un
espace de Fr\'echet. Elle ne d\'epend ni de l'atlas $\{ (U,\varphi
)\} $ ni de la famille ${\cal C}$ ni de la suite croissante de
compacts $K_n$. C'est la topologie $C^{0,\infty }$ sur $A(M)=C^{0,\infty }(M)$.
\smallskip
Pour cette topologie, le sous-espace ${\cal H}_{\cal F}(M)$ des
fonctions ${\cal F}$-holomorphes est ferm\'e. Notons que, pour
d\'efinir la topologie induite sur   ${\cal H}_{\cal F}(M)$, on peut
faire l'\'economie des d\'eriv\'ees le long des feuilles : il suffit
de consid\'erer l'op\'erateur diff\'erentiel $D^0$.
\smallskip

La topologie $C^{0,\infty }$ se d\'efinit de mani\`ere analogue sur
les espaces $A^{pq}({\cal F})$, vu que tout \'el\'ement  $\alpha
\in \Omega^{pq}({\cal F})$ s'\'ecrit, sur une carte locale
$(U,\varphi )$, sous la forme $\alpha =fd\overline z$, $\alpha =fdz$
ou $\alpha =fdz\wedge d\overline z$ avec $f\in A(U)$.

\smallskip
\noindent {\tfpp 3.2. Z\'eros et p\^oles}
\smallskip
On travaillera sur un ouvert $U$ de $\Bbb C\times {\Bbb R}$ muni de
son feuilletage complexe canonique ${\cal F}$. (Le facteur $\Bbb R$
pourrait \^etre remplac\'e par n'importe quelle vari\'et\'e
diff\'erentiable et notamment par ${\Bbb R}^n$.)
\medskip
Soient $f:U\longrightarrow \Bbb C$ une fonction ${\cal
F}$-holomorphe et $Z$ l'ensemble de ses z\'eros. La restriction de
$f$ \`a chaque feuille $F$  est une fonction holomorphe ; par suite,
si $f:F\longrightarrow \Bbb C$ n'est pas identiquement nulle, par le
principe des z\'eros isol\'es,  $Z\cap F$ est une partie discr\`ete
de $F$.  Donc en un point de $Z\cap F$ o\`u $f$ n'est pas
idendiquement nulle, $Z\cap F$ est ``transverse" \`a $F$.

Une fonction $f:U\longrightarrow \Bbb C$ est dite ${\cal F}$-{\it
m\'eromorphe},  si sa restriction \`a chaque feuille est une
fonction m\'eromorphe. Notons ${\cal P}$ l'ensemble des p\^oles de
$f$ ; alors, comme pour les z\'eros, l'intersection de ${\cal P}$
avec toute feuille $F$ est un ensemble discret de $F$.

Une fonction ${\cal F}$-holomorphe $f:U\longrightarrow \Bbb C$
(resp. ${\cal F}$-m\'eromorphe) \'etant simplement continue sur
$U$ (resp. sur $U\setminus {\cal P}$), on ne peut malheureusement
pas dire plus ni sur  l'ensemble de ses z\'eros ni celui de ses
  p\^oles. Nous ferons simplement des remarques lorsque $Z$ et
${\cal P}$ poss\`edent une structure de vari\'et\'e $C^\infty $ ; de
telles fonctions existent bien s\^ur : si $\varphi  :]-\eta ,\eta
[\longrightarrow U\subset \Bbb C$ est une courbe diff\'erentiable,
alors les fonctions $f(z,t)=z-\varphi  (t)$ et $g(t)={{1}\over
{f(z,t)}}$ sont ${\cal F}$-holomorphes et ont $\{ (z,t)\in {\Bbb
C}\times {\Bbb R}: z=\varphi (t)\hbox{ avec } t\in ]-\eta ,\eta [\}
$  respectivement comme ensemble de z\'eros et ensemble de p\^oles.
Ceci permet de d\'efinir localement des fonctions du m\^eme type sur
des feuilletages holomorphes $(M,{\cal F})$.
\smallskip
Soit maintenant $\Sigma $ une petite sous-vari\'et\'e  transverse
\`a ${\cal F}$ ; elle  peut \^etre consid\'er\'ee comme  le
graphe d'une application $z_0 :t\in ]-\eta ,\eta[\longmapsto
z_0(t)\in \Bbb C$ de classe $C^\infty $. Soit $V$ un voisinage
ouvert relativement compact de $\Sigma $ dont chaque section $V^t$
est un disque centr\'e en $z_0(t)$. Alors, sur $V\setminus \Sigma $,
la fonction  $f$ admet un d\'eveloppement de Laurent : $$f(z,t)=
\sum_{n=0}^\infty a_n(t)(z-z_0(t))^n+\sum_{m=1}^\infty
{{b_m(t)}\over {(z-z_0(t))^m}}$$  o\`u les coefficients $a_n$ et
$b_m$ sont donn\'es, comme dans le cas classique, par les formules
int\'egrales :
$$a_n(t)={{1}\over {2i\pi }}\int_{\gamma_1^t}{{f(\xi
,t)}\over {(\xi -z_0(t))^{n+1}}}d\xi \;\;\hbox{et}\;\;
b_m(t)={{1}\over {2i\pi }}\int_{\gamma_2^t}(\xi -z_0(t))^{n-1}f(\xi
,t)d\xi . \leqno{(7)}$$ $\gamma_1^t$ et $\gamma_2^t$ sont respectivement le
grand cercle et le petit cercle d'une couronne contenant le point
$(z,t)$ dans la section $V^t$ de $V$. Le point $(z_0,t_0)\in \Sigma
$ est une singularit\'e si l'un au moins des $b_m(t_0)$ est non nul
; s'il existe $m_0\geq 1$ tel que $b_{m_0}(t_0)\neq 0$ et
$b_m(t_0)=0$ pour $m>m_0$, on dira que $(z_0,t_0)\in \Sigma $ est un
{\it p\^ole} de $f$ d'{\it ordre} $m_0$ ; s'il existe une infinit\'e
de  $b_m(t_0)$ non nuls, on dira que $(z_0,t_0)\in \Sigma $ est une
{\it singularit\'e essentielle} ; si $a_0(t_0)$ et tous les
$b_m(t_0)$ sont nuls, on dira que $(z_0,t_0)$ est un {\it z\'ero} de
$f$ ; sa {\it multiplicit\'e} est par d\'efinition le plus petit
entier $n\geq 1$ tel que $a_n(t_0)\neq 0$.  Comme les $a_n$ et les
$b_m$ sont des fonctions continues en $t$, si $(z_0,t_0)\in \Sigma
$ est un point singulier de $f$ ({\it i.e.} un p\^ole ou une
singularit\'e essentielle) par continuit\'e, il existe $\delta >0$
tel que pour $\vert t-t_0\vert <\delta $, le point $(z_0(t),t)$ est
aussi singulier. L'ensemble singulier de $f$ est donc   une
transversale \`a ${\cal F}$ qui est ouverte.

Ces remarques nous permettent de montrer facilement la proposition
suivante dont on fera usage dans le calcul explicite de l'exemple 2.

\smallskip
\noindent {\tfpp 3.3. Proposition.}  {\it Soit $f:M\longrightarrow
\Bbb C$ une fonction continue et
${\cal F}$-holomorphe en dehors d'une r\'eunion discr\`ete de
sous-vari\'et\'es transverses $\Sigma_j $. Alors} :

\item{i)} {\it chaque $\Sigma_j $ est  ouverte} ;

\item{ii)} {\it si chacune des  $\Sigma_j $ est r\'eduite \`a un   point, $f$ se
prolonge en une fonction ${\cal F}$-holomorphe sur toute la vari\'et\'e
$M$. (C'est un phénoménome type Hartogs.)}
\medskip
{\tfp 4. Le feuilletage complexe affine de Reeb}
\medskip
\smallskip
\noindent {\tfpp 4.1. Construction}
\smallskip
$\bullet $  On pose $\widetilde M={\Bbb C}\times {\Bbb R}^n\setminus \{ (0,0)\} $ et on consid\`ere le feuilletage
$\widetilde{{\cal F}}$ d\'efini par le syst\`eme diff\'erentiel
$dt_1=\ldots =dt_n=0$
o\`u $(z,t_1,\ldots ,t_n)$ sont les coordonn\'ees d'un
point $(z,t)$ dans $\widetilde M$. Les feuilles de $\widetilde{\cal F}$
sont holomorphiquement \'equivalentes \`a ${\Bbb C}$ sauf
celle qui correspond \`a $t=0$ qui est ${\Bbb C}^\ast $.  Soit
$\phi :\widetilde M \longrightarrow \widetilde M$ le
diff\'eomorphisme de $M$ d\'efini par $\phi (z,t) = (\lambda z,\lambda
t)$ o\`u
$\lambda \in ]0,1[$. L'action de ${\Bbb Z}$ engendr\'ee par $\phi $
est libre, propre et discontinue ; le quotient $M=\widetilde M/\phi$
est diff\'eomorphe  \`a la {\it vari\'et\'e de Hopf} r\'eelle
${\Bbb S}^1\times {\Bbb S}^{n+1}$.   Le feuilletage
$\widetilde{\cal F}$ est invariant par $\phi $ et induit  un
feuilletage complexe ${\cal F}_{\lambda }$ de dimension $1$
sur $M$. Les feuilles sont des copies de ${\Bbb C}$ sauf celle qui
correspond \`a $t=0$ qui est une  {\it
courbe elliptique} $C_{\lambda }$ dont la structure complexe est
donn\'ee par celle de la couronne $\{ z\in {\Bbb C}\; : \; \vert
\lambda \vert <\vert z\vert <1\} .$ Tout isomorphisme de feuilletages
$f:(M,{\cal F}_{\lambda })\longrightarrow (M,{\cal F}_{\lambda' })$
induit un biholomorphisme de $C_\lambda $ sur
$C_{\lambda'}$. Donc si $\lambda  \neq \lambda'$,
${\cal F}_{\lambda }$ n'est pas isomorphe \`a ${\cal
F}_{\lambda'}$.
\smallskip
$\bullet $  Cherchons le groupe $G({\cal F}_\lambda)$ des automorphismes de ${\cal F}_\lambda$
sur $M=\widetilde M/\phi$ (dans le cas  o\`u la variété  $\widetilde M$ est $\Bbb
C\times \Bbb  R$). Un \'el\'ement $f$ de $G({\cal F}_\lambda)$ est donn\'e par
un automorphisme $\widetilde f:\widetilde M\longrightarrow
\widetilde M$ de ${\cal F}$ commutant \`a l'action de $\phi $ ; il
s'\'ecrit $\widetilde f(z,t)= (f_1(z,t),f_2(t))$ o\`u $f_1$ est
holomorphe en $z$ et commute \`a la multiplication $z\longmapsto
\lambda z$ et $f_2$ est un diff\'eomorphisme de $\Bbb  R$ ; $f_1$ est
n\'ecessairement de la forme $f_1(z,t)=a(t)z$ o\`u $a(t)\in
\hbox{GL}(1,\Bbb  C)=\Bbb C^\ast $ d\'ependant diff\'erentiablement
de $t$. D'autre part,  comme $\Bbb  C^\ast $ n'est pas \'equivalent
\`a $\Bbb  C$, $f_2$ doit fixer $0$ et commuter \`a l'homoth\'etie
$\phi_2:t\longmapsto \lambda t$ {\it i.e.} $f_2(\lambda t) = \lambda
f_2(t)$ ; il est alors du type $f_2(t) = bt$ o\`u $b\in \Bbb  R^\ast
$. Le groupe $G({\cal F}_\lambda )$ est donc celui
  des transformations  de la forme
$(z,t)\longmapsto (a(t)z,bt)$ o\`u $a\in C^\infty (\Bbb  R,\Bbb
C^\ast )$ et $b\in \Bbb  R^\ast $.

\smallskip
La vari\'et\'e ici est $M={\Bbb S}^1\times {\Bbb S}^{n+1}$ qui est le quotient de
$\widetilde M={\Bbb C}\times {\Bbb R}^n\setminus \{ (0,0)\} $ par
l'action de ${\Bbb Z}$ engendr\'ee par l'automorphisme
$(z,t)\in \widetilde M\buildrel \gamma \over \longrightarrow (\lambda z,\lambda t)\in \widetilde M.$
Elle sera munie du feuilletage ${\cal F}_{\lambda }$ d\'efini qu'on notera simplement ${\cal F}$.
\smallskip
\noindent {\tfpp 4.2. Th\'eor\`eme principal.} {\it Les espaces vectoriels  $H_{\cal F}^{00}(M)={\cal H}_{\cal F}(M)$ et
 $H_{\cal F}^{01}(M)$ sont isomorphes \`a la droite complexe ${\Bbb C}$. }
\smallskip
$\bullet $ R\'eglons d'abord le cas de ${\cal H}_{\cal F}(M)$. Un \'el\'ement de l'espace ${\cal H}_{\cal F}(M)$
est une fonction ${\cal F}$-holomorphe $f$ sur $M$, donc une fonction sur $\widetilde M$ telle que $f(\lambda z,\lambda t)=f(z,t)$
et qui est ${\cal F}$-holomorphe. Mais, d'apr\`es la Proposition 3.3, $f$ s'\'etend en une fonction ${\cal F}$-holomorphe
sur ${\Bbb C}\times {\Bbb R}^n$. D'autre part :
$$f(z,t)=\sum_{\ell =0}^\infty f_\ell (t)z^\ell $$
o\`u les $f_\ell $ sont des fonctions continues de $t\in {\Bbb R}^n$ (donn\'ees par les formules int\'egales de gauche dans (7)). Comme $f$ v\'erifie $f(\lambda z,\lambda t)=f(z,t)$,
les $f_\ell $ doivent satisfaire la relation $\lambda^{-\ell }f_\ell (t)=f_\ell (\lambda t)$. Mais, pour $\ell \neq 0$,  cette
relation force $\vert f_\ell \vert $ \`a tendre vers $+\infty $ quand $\vert t\vert \to 0$ et ne saurait donc \^etre continue
pour $\ell \neq 0$. Donc $f$ est r\'eduite \`a la fonction $f_0$ ; celle-ci v\'erifiant $f_0(\lambda t)=f_0(t)$ doit \^etre en fait la constante
$f_0(0)$.
\smallskip
Pour  $H_{\cal F}^{01}(M)$, la d\'emonstration se fera pour $n=1$. (Elle est exactement la m\^eme pour $n\geq 2$.) Nous allons commencer
par d\'ecrire les espaces $A_{\cal F}^{pq}(M)$ des formes feuillet\'ees de type $(p,q)$.
\smallskip
$\bullet $ Une forme feuillet\'ee de type $(0,1)$ sur $\widetilde M$ (resp. de type $(1,0)$) s'\'ecrit $\alpha =f(z,t)d\overline z$ (resp. $\beta =g(z,t)dz$)
o\`u $f$ et $g$ sont des \'el\'ements de $A(M)$. L'action de $\gamma $ sur $\alpha $ (resp. $\beta $) est donn\'ee par $\gamma^\ast (\alpha ) =\lambda f(\lambda z, \lambda t)d\overline z$
(resp.  $\gamma^\ast (\beta ) =\lambda g(\lambda z, \lambda t)d\overline z$). Les deux formes $\alpha $ et $\beta $ sont donc invariantes par $\gamma $ si les
fonctions $f$ et $g$ v\'erifient les relation fonctionnelles :
$$\lambda f(\lambda z,\lambda t)=f(z,t)\hskip0.5cm \hbox{et} \hskip0.5cm \lambda g(\lambda z,\lambda t)=g(z,t).\leqno{(8)}$$
\smallskip
$\bullet $ Une forme feuillet\'ee de type $(1,1)$ sur $\widetilde M$ s'\'ecrit $\eta =h(z,t)dz\wedge d\overline z$. La condition d'invariance par $\gamma $ impose \`a
la fonction $h$ de v\'erifier cette fois-ci $\lambda^2h(\lambda z,\lambda t)=h(z,t).$

$\bullet $ Donnons explicitement quelques exemples de ces formes feuillet\'ees, par exemple de type $(0,1)$, situation
qui va le plus nous int\'eresser par la suite. On doit donc trouver une fonction $f\in A(M)$ telle que $\lambda f(\lambda z,\lambda t)=f(z,t)$.
Il est \'evident que :
$$a(z,t)={1\over {\sqrt{z\overline z+t^2}}}$$
en est une. Si $f$ est une autre fonction (tout \`a fait quelconque) et v\'erifiant (8),
la fonction $u={f\over {f_0}}$
(bien d\'efinie car $f_0$ ne s'annule nulle part) v\'erifie $u(\lambda z,\lambda t)=u(z,t)$, donc une fonction sur la vari\'et\'e compacte
$M=\widetilde M/\langle \gamma \rangle $ et par suite born\'ee. Toute $(0,1)$-forme feuillet\'ee sur $M$ s'\'ecrit
donc $\alpha (z,t)=a(z,t)u(z,t)d\overline z$ o\`u $u$ est une fonction sur $\widetilde M$ invariante par l'action du diff\'eomorphisme
$\gamma $.

\smallskip
\noindent {\tfpp 4.3. D\'emonstration du th\'eor\`eme principal}
\smallskip
\underbar{Premi\`ere \'etape}
\smallskip
\noindent Soient $R$ et $\varepsilon $ deux r\'eels strictement positifs tels que $R+\varepsilon <\lambda^{-1}R$ ; on note  $L$ le disque ferm\'e de ${\Bbb C}$ de rayon $R$
et $\Omega $ un $\varepsilon $-voisinage ouvert de $L$. Pour tout $j\in {\Bbb N}$, on pose $K_j=\lambda^{-j}L_0$ et $\Omega_j=\lambda^{-j}\Omega_0$
o\`u $L_0=L\times {\Bbb R} $ et $\Omega_0=\Omega \times {\Bbb R} .$
\smallskip

Soit $\rho_0:z\in {\Bbb C}\longmapsto \rho_0(z)\in {\Bbb R}_+$ une fonction $C^\infty $ \`a support compact, ne d\'ependant que de $\vert z\vert $
 et identiquement \'egale \`a $1$ sur $\Omega $. Pour
tout $j\in {\Bbb N}$, la fonction $\rho_j(\xi )=\rho_0(\lambda^j\xi )$ est $C^\infty $, \`a support compact et vaut $1$ identiquement
sur l'ouvert $\lambda^{-j}\Omega $.
 Soit $\phi_j:{\Bbb C}\times {\Bbb R}\longrightarrow {\Bbb R}_+$ la fonction d\'efinie par $\phi_j(\xi ,t)=\rho_j(\xi )$ et :
$$\psi_j(\xi ,t)=\cases{\phi_j(\xi ,t)&\hbox{si $j=0$}\cr \phi_j(\xi ,t)-\phi_{j-1}(\xi ,t) &\hbox{si $j\geq 1$.}}$$
Comme  $\phi_j(\xi ,t)$ et $\psi_j(\xi ,t)$ ne d\'ependent pas de $t$, on les notera simplement $\phi_j(\xi )$ et $\psi_j(\xi )$.
Les relations suivantes sont bien s\^ur imm\'ediates mais comme elles nous seront tr\`es utiles, nous les rappelons et les mettrons
bien en vue :
$$\cases{\phi_j(\xi )\hskip0.2cm =\phi_0(\lambda^j\xi )\hskip0.2cm \hbox{(par d\'efinition de $\phi_j$)}\cr \psi_j(\lambda \xi )=\psi_{j+1}(\xi ).}\leqno{(9)}$$
\smallskip
Pour tout $j\in {\Bbb N}$ et tout $t\in {\Bbb R}$, la section $K_j^t$ de $K_j$ est un
${\cal F}$-compact (son intersection avec toute feuille est un compact)
contenu dans l'int\'erieur de $K_{j+1}$. On a :
$$\bigcup_{j=0}^\infty K_j={\Bbb C}\times {\Bbb R} \hskip0.5cm \hbox{et} \hskip0.5cm  \sum_{j=0}^\infty \psi_j=1.$$
\smallskip
\underbar{Deuxi\`eme \'etape}
\smallskip
\noindent On se donne une $(0,1)$-forme feuillet\'ee $\alpha =f(z,t)d\overline z$ o\`u $f$ v\'erifie la relation fonctionnelle (8).
Pour une raison \'evidente de degr\'e, cette forme est $\overline \partial_{\cal F}$-ferm\'ee. D'autre part, comme
$\alpha $ s'\'ecrit $\alpha =a(z,t)u(z,t)d\overline z$ avec $u$ invariante par $\gamma $, la fonction $f$ a la croissance
de la fonction $a$ qui est localement int\'egrable.
Pour tout $(z,t)\in {\Bbb C}\times {\Bbb R}\setminus \{ (0,0)\} =\widetilde M$, la quantit\'e qui suit existe :
$$h_0(z,t)={1\over {2i\pi }}\int_{\Bbb C}{{\psi_0(\xi )f(\xi ,t)}\over {\xi -z}}d\xi \wedge d\overline \xi .$$
La fonction $h_0$ est donc bien d\'efinie, continue en $(z,t)$ et $C^\infty $ en $z$ ({\it cf.} [H\"o] Theorem 1.2.2 page 3). En utilisant la formule
int\'egrale de Cauchy, on montre facilement que $h_0$ v\'erifie l'\'equation $\overline \partial_{\cal F}h_0=\psi_0f$
sur $\widetilde M$.
\smallskip
Soit $j\in {\Bbb N}^\ast $. Le support de $\psi_j$ ne contient pas l'origine $(0,0)$. Comme pr\'ec\'edemment, on pose
pour tout $(z,t)\in \widetilde M$ :
$$h_j(z,t)={1\over {2i\pi }}\int_{\Bbb C}{{\psi_j(\xi )f(\xi ,t)}\over {\xi -z}}d\xi \wedge d\overline \xi .$$
Comme pour $h_0$, la fonction $h_j$ est bien d\'efinie, continue en $(z,t)$ et $C^\infty $ en $z$.
Elle   v\'erifie l'\'equation $\overline \partial_{\cal F}h_j=\psi_jf$.

\smallskip
\underbar{Troisi\`eme \'etape}
\smallskip

\noindent Il nous reste \`a recoller toutes les solutions partielles que nous avons obtenues.
Comme $\psi_j=0$ sur $\Omega_{j-1}$,
$h_j$ y est ${\cal F}$-holomorphe ({\it cf.} [H\"o]). Les sections $\Omega_{j-1}^t$ \'etant des disques ouverts,
on peut d\'evelopper $h_j$ en s\'erie enti\`ere :
$$h_j(z,t)=\sum_{n=0}^\infty a_n(t)z^n\leqno{(10)}$$
qui converge pour la m\'etrique $\delta $ sur l'ouvert  $\Omega_{j-1}$.
En tronquant de fa\c con ad\'equate la s\'erie (10), on obtient une fonction $v_j$, ${\cal F}$-holomorphe
sur ${\Bbb C}\times {\Bbb R}$ (c'est un polyn\^ome en $z$) et telle que :
$$\delta (h_j,v_j)<{1\over {2^j}}. \leqno{(11)}$$
Soit  $\widetilde h:\widetilde M\longrightarrow {\Bbb C}$   la fonction d\'efinie   par :
$$\widetilde h(z,t)=h_0(z,t)+\sum_{j=1}^\infty \left( h_j(z,t)-v_j(z,t)\right) .$$
En vertu de l'in\'egalit\'e  (11), la s\'erie converge uniform\'ement au sens de la m\'etrique $\delta $ ; la fonction $\widetilde h$ est donc
continue en $(z,t)$ et de classe $C^\infty $ en $z$. En plus, comme l'op\'erateur $\overline \partial_{\cal F}$
est continu pour la topologie $C^{0,\infty }$, on a :
$$\eqalign{\overline \partial_{\cal F}\widetilde h&=\overline \partial_{\cal F}\left( h_0(z,t)+\sum_{j=1}^\infty \left( h_j(z,t)-v_j(z,t)\right) \right) \cr
& = \overline \partial_{\cal F}h_0(z,t)+\sum_{j=1}^\infty \left( \overline \partial_{\cal F}h_j(z,t)-\overline \partial_{\cal F}v_j(z,t)\right)  \cr
& = \sum_{j=0}^\infty \overline \partial_{\cal F}h_j(z,t)\cr
&= \sum_{j=0}^\infty \psi_j(z,t)f(z,t)\cr
&= f}$$
qui montre bien que $h$ est une solution de l'\'equation $\overline \partial_{\cal F}\widetilde h=f$. Mais, a priori, elle ne d\'efinit pas
une solution au probl\`eme sur la vari\'et\'e quotient $M=\widetilde M/\langle \gamma \rangle $ : elle ne v\'erifie pas forc\'ement la
condition d'invariance $\widetilde h(\lambda z,\lambda t)=\widetilde h(z,t)$. Pour en obtenir une, on corrige $\widetilde h$ en lui rajoutant
une fonction ${\cal F}$-holomorphe $K(z,t)$ de telle sorte que $h(z,t)=  \widetilde h(z,t)+K(z,t)$, qui v\'erifie encore l'\'equation
$\overline \partial_{\cal F} h=f$, soit $\gamma $-invariante {\it i.e.}  $h(\lambda z,\lambda t)=h(z,t)$, ce qui impose \`a $K$ de v\'erifier
l'{\it \'equation cohomologique}  :
$$K(z,t)-K(\lambda z,\lambda t)=H(z,t)\leqno{(12)}$$
o\`u $H(z,t)= \widetilde h(\lambda z,\lambda t)-\widetilde h(z,t)$. Nous avons donc \`a r\'esoudre l'\'equation (12)
o\`u l'inconnue est la fonction ${\cal F}$-holomorphe $K$. Remarquons que $H$ est continue et  ${\cal F}$-holomorphe sur $\widetilde M$ ; en effet :
$$\eqalign{\overline \partial_{\cal F}H(z,t)&=\overline \partial_{\cal F}\left( \widetilde h(\lambda z,\lambda t)-\widetilde h(z,t)\right) \cr
&=\lambda(\overline \partial_{\cal F}\widetilde h)(\lambda z,\lambda t)-\overline \partial_{\cal F}\widetilde h(z,t) \cr
&=\lambda f(\lambda z,\lambda t)-f(z,t)\cr
&= f(z,t)-f(z,t)\cr
&=0.}$$
\'Etant ${\cal F}$-holomorphe sur l'ouvert $\widetilde M$, $H$ l'est sur l'espace ${\Bbb C}\times {\Bbb R}$ tout entier en vertu de la Proposition 3.3.
\smallskip
\underbar{Quatri\`eme \'etape}
\smallskip
\noindent Formellement la fonction $K(z,t)=\displaystyle \sum_{n=0}^\infty H(\lambda^nz,\lambda^nt)$ est solution de l'\'equation (12). Il ne reste donc plus
qu'\`a montrer que cette  s\'erie converge pour la topologie $C^{0,\infty }$ pour qu'elle d\'efinisse effectivement une fonction  continue
et ${\cal F}$-holomorphe.
\smallskip
Une condition n\'ecessaire de l'existence de $K$  est $H(0,0)=0$. Avant d'examiner comment elle est v\'erifi\'ee,
explicitons la quantit\'e   :
$$\eqalign{H(z,t)&=\widetilde h(\lambda z,\lambda t)-\widetilde h(z,t)\cr
&=h_0(\lambda z,\lambda t)-h_0(z,t)
  + \sum_{j=1}^\infty \left\{  h_j(\lambda z,\lambda t)-h_j(z,t)\right\}
  - \sum_{j=1}^\infty \left\{  v_j(\lambda z,\lambda t)-v_j(z,t)\right\} .}$$
\smallskip
\noindent $\bullet $ Commen\c cons par $h_0(\lambda z,\lambda t)-h_0(z,t)$. On a :
$$\eqalign{h_0(\lambda z,\lambda t)-h_0(z,t)&={1\over {2i\pi }}\int_{\Bbb C}{{\psi_0(\xi )f(\xi ,\lambda t)}\over {\xi -\lambda z}}d\xi \wedge d\overline \xi -
{1\over {2i\pi }}\int_{\Bbb C}{{\psi_0(\xi )f(\xi ,t)}\over {\xi -z}}d\xi \wedge d\overline \xi \cr
&= {1\over {2i\pi }}\int_{\Bbb C}{{\psi_0(\lambda \xi )f(\xi ,t)}\over {\xi -z}}d\xi \wedge d\overline \xi -
{1\over {2i\pi }}\int_{\Bbb C}{{\psi_0(\xi )f(\xi ,t)}\over {\xi -z}}d\xi \wedge d\overline \xi \cr
&={1\over {2i\pi }}\int_{\Bbb C}{{\psi_1(\xi )f(\xi , t)}\over {\xi -z}}d\xi \wedge d\overline \xi .}\leqno{(13)}$$
Le passage de la premi\`ere ligne \`a la deuxi\`eme se fait par chagement de variable $\xi \longmapsto \lambda \xi $ et utilise la relation $\lambda f(\lambda \xi ,\lambda t)=
f(\xi ,t)$ et le passage de la deuxi\`eme ligne \`a la troisi\`eme les relations (9).
\smallskip
\noindent $\bullet $ Un calcul similaire donne :
$$h_j(\lambda z,\lambda t)-h_j(z,t)={1\over {2i\pi }}\int_{\Bbb C}{{(\psi_j(\lambda \xi )-\psi_j(\xi ))
f(\xi ,  t)}\over {\xi -z}}d\xi \wedge d\overline \xi .$$
En utilisant la deuxi\`eme relation de (9), on obtient :
$$h_j(\lambda z,\lambda t)-h_j(z,t)={1\over {2i\pi }}\int_{\Bbb C}{{(\psi_{j+1}(\xi )-\psi_j(\xi ))
f(\xi ,  t)}\over {\xi -z}}d\xi \wedge d\overline \xi .\leqno{(14)}$$
\smallskip
\noindent $\bullet $ Les relations (13) et (14) donnent finalement :
$$\sum_{j=0}^N\left\{ h_j(\lambda z,\lambda t)-h_j(z,t)\right\} =
{1\over {2i\pi }}\int_{\Bbb C}{{\psi_N(\xi )
f(\xi ,t)}\over {\xi -z}}d\xi \wedge d\overline \xi .$$
Comme l'\'evaluation de la quantit\'e $\displaystyle \sum_{j=1}^\infty \left\{  v_j(\lambda z,\lambda t)-v_j(z,t)\right\} $
en $(0,0)$ est nulle (pour tout $j\geq 1$, la fonction $v_j$ est d\'efinie en $(0,0)$), on obtient :
$$H(0,0)=\lim_{N\to +\infty }\left( {1\over {2i\pi }}\int_{\Bbb C}{{\psi_N(\xi )
f(\xi ,0)}\over {\xi }}d\xi \wedge d\overline \xi \right) .$$
\smallskip
\noindent $\bullet $ Pour finir cette \'etape, montrons que la suite de nombres complexes $(I_N)_{N\geq 1}$ o\`u $I_N$ est
donn\'e par :
$$I_N={1\over {2i\pi }}\int_{\Bbb C}{{\psi_N(\xi )
f(\xi ,0)}\over {\xi }}d\xi \wedge d\overline \xi $$
est constante. Ceci r\'esulte du calcul imm\'ediat qui suit, qui utilise la deuxi\`eme des relations (9)
et l'invariance $\lambda f(\lambda \xi ,0)=f(\xi ,0)$ de $f$ (pour $\xi \neq 0$ bien s\^ur). En effet :
$$\eqalign{I_{N+1}&= {1\over {2i\pi }}\int_{\Bbb C}{{\psi_{N+1}(\xi )
f(\xi ,0)}\over {\xi }}d\xi \wedge d\overline \xi \cr
&={1\over {2i\pi }}\int_{\Bbb C}{{\psi_N(\lambda \xi )
f(\xi ,0)}\over {\xi }}d\xi \wedge d\overline \xi \cr
&= {1\over {2i\pi }}\int_{\Bbb C}{{\psi_N(\zeta )
f({\zeta \over \lambda } ,0)}\over {{\zeta \over \lambda } }}{{d\zeta  \wedge d\overline \zeta }\over {\lambda^2}}\cr
&={1\over {2i\pi }}\int_{\Bbb C}{{\psi_N(\zeta )
f(\zeta ,0)}\over {\xi }}d\zeta \wedge d\overline \zeta \cr
&=I_N.}$$
Par suite, la condition $H(0,0)=0$ est \'equivalente \`a :
$${1\over {2i\pi }}\int_{\Bbb C}{{\psi_1(\xi )
f(\xi ,0)}\over {\xi }}d\xi \wedge d\overline \xi =0.\leqno{(15)}$$
\smallskip
\underbar{Cinqui\`eme \'etape}
\smallskip
\noindent Nous avons vu que la condition $H(0,0)=0$ est n\'ecessaire \`a l'existence de la fonction $K$. Montrons maintenant qu'elle est
aussi suffisante. Ce sera le cas si on montre la convergence de la s\'erie qui suit pour la topologie $C^{0,\infty }$ :
$$\eqalign{\sum_{n=0}^\infty H\left( \lambda^nz,\lambda^nt\right) &=\sum_{n=0}^\infty \left( h_0(\lambda^{n+1}z,\lambda^{n+1}t)-h_0(\lambda^nz,\lambda^nt)\right) \cr
& \hskip0.4cm + \sum_{n=0}^\infty \sum_{j=1}^\infty \left\{ h_j(\lambda^{n+1}z,\lambda^{n+1}t)-h_j(\lambda^nz,\lambda^nt)\right\}  \cr
& \hskip0.4cm -\sum_{n=0}^\infty \sum_{j=1}^\infty \left\{  v_j(\lambda^{n+1}z,\lambda^{n+1}t)-v_j(\lambda^nz,\lambda^nt)\right\} .}\leqno{(16)}$$
Un calcul simple, utilisant le fait que $f$ v\'erifie la relation fonctionnelle $\lambda f(\lambda \xi ,\lambda t)=f(\xi ,t)$, montre que :
$$h_0(\lambda^{n+1}z,\lambda^{n+1}t)-h_0(\lambda^nz,\lambda^nt)={1\over {2i\pi }}\int_{\Bbb C}
{{\left\{ \psi_0(\lambda^{n+1}\xi )-\psi_0(\lambda^n\xi )\right\} f(\xi ,t)}\over {\xi -z}}d\xi \wedge d\overline \xi $$
et :
$$h_j(\lambda^{n+1}z,\lambda^{n+1}t)-h_j(\lambda^nz,\lambda^nt)={1\over {2i\pi }}\int_{\Bbb C}
{{\left\{ \psi_j(\lambda^{n+1}\xi )-\psi_j(\lambda^n\xi )\right\} f(\xi ,t)}\over {\xi -z}}d\xi \wedge d\overline \xi .$$
Pour montrer la convergence de la s\'erie $K(z,t)=\displaystyle \sum_{n=0}^\infty H(\lambda^nz,\lambda^nt)$, on va expliciter
et simplifier l'expression des trois s\'eries qui composent le membre de droite de la relation (16). Formellement on a :
$$\sum_{n=0}^\infty \left( h_0(\lambda^{n+1}z,\lambda^{n+1}t)-h_0(\lambda^nz,\lambda^nt)\right) =\int_{\Bbb C}
{{\sum_{n=0}^\infty \left\{ \psi_0(\lambda^{n+1}\xi )-\psi_0(\lambda^n\xi )\right\} f(\xi ,t)}\over {2i\pi(\xi -z)}}d\xi \wedge d\overline \xi .$$
Soit $N$ un entier naturel positif. Comme  :
$$\sum_{n=0}^N\left\{ \psi_0(\lambda^{n+1}\xi )-\psi_0(\lambda^n\xi )\right\} =\psi_0(\lambda^{N+1}\xi )-\psi_0(\xi )\leqno{(17)}$$
on a :
$$\sum_{n=0}^N\left\{ h_0(\lambda^{n+1}z,\lambda^{n+1}t)-h_0(\lambda^nz,\lambda^nt)\right\}  =
{1\over {2i\pi }}\int_{\Bbb C}{{ \left( \psi_0(\lambda^{N+1}\xi )-\psi_0(\xi )\right) f(\xi ,t)}\over {\xi -z}}d\xi \wedge d\overline \xi .$$
De la m\^eme
mani\`ere, nous allons nous occuper de la s\'erie double  :
$$\displaystyle \sum_{j=1}^\infty \left\{  h_j(\lambda^{n+1}z,\lambda^{n+1}t)-h_j(\lambda^nz,\lambda^nt)\right\} .$$
D'abord on a :
$$\eqalign{\psi_j(\lambda^{n+1}\xi )&=\phi_j(\lambda^{n+1}\xi )-\phi_{j-1}(\lambda^{n+1}\xi )\cr
&=\phi_0(\lambda^{n+j+1}\xi )-\phi_0(\lambda^{n+j}\xi ).}$$
En sommant sur $n\in {\Bbb N}$ de $0$ \`a $N$, on obtient :
$$\sum_{n=0}^N\psi_j(\lambda^{n+1}\xi )=\phi_0(\lambda^{j+N}\xi )-\phi_0(\lambda^j\xi ).$$
De fa\c con similaire, on \'etablit l'\'egalit\'e :
$$\sum_{n=0}^N\psi_j(\lambda^n\xi )=\phi_0(\lambda^{j+N}\xi )-\phi_0(\lambda^{j-1}\xi ).$$
Et par suite :
$$\eqalign{\sum_{n=0}^N\left\{ \psi_j(\lambda^{n+1}\xi )-\psi_j(\lambda^n\xi )\right\} &=-\phi_0(\lambda^j\xi )+\phi_0(\lambda^{j-1}\xi )\cr
&= -\phi_j(\xi )+\phi_{j-1}(\xi )\cr
&= -\psi_j(\xi ).}$$
Finalement :
$$\sum_{n=0}^N\left\{ h_j(\lambda^{n+1}z,\lambda^{n+1}t)-h_j(\lambda^nz,\lambda^nt)\right\}  =
{1\over {2i\pi }}\int_{\Bbb C}{{(-\psi_j(\xi )) f(\xi ,t)}\over {\xi -z}}d\xi \wedge d\overline \xi .$$
Apr\`es sommation sur $j\in {\Bbb N}^\ast $, on obtient :
$$\sum_{j=1}^\infty \sum_{n=0}^N\left\{ h_j(\lambda^{n+1}z,\lambda^{n+1}t)-h_j(\lambda^nz,\lambda^nt)\right\}  =
{1\over {2i\pi }}\int_{\Bbb C}{{(\psi_0(\xi)-1) f(\xi ,t)}\over {\xi -z}}d\xi \wedge d\overline \xi .$$
D'o\`u :
$$\sum_{n=0}^N\left\{ h_0(\lambda^{n+1}z,\lambda^{n+1}t)-h_0(\lambda^nz,\lambda^nt)\right\} +
\sum_{j=1}^\infty \sum_{n=0}^N\left\{ h_j(\lambda^{n+1}z,\lambda^{n+1}t)-h_j(\lambda^nz,\lambda^nt)\right\} =$$
$${1\over {2i\pi }}\int_{\Bbb C}{{\left( (\psi_0(\xi )-1)+ (\psi_0(\lambda^{N+1}\xi )-\psi_0(\xi ))\right) f(\xi ,t)}\over {\xi -z}}d\xi \wedge d\overline \xi .$$
On fait tendre $N$ vers $+\infty $ ; le terme $\psi_0(\lambda^{N+1}\xi )$ tend vers $\psi_0(0)=1$ (puisque $0<\lambda <1$). Par suite :
$$\sum_{n=0}^\infty \left\{ h_0(\lambda^{n+1}z,\lambda^{n+1}t)-h_0(\lambda^nz,\lambda^nt)\right\} +
\sum_{j=1}^\infty \sum_{n=0}^\infty \left\{ h_j(\lambda^{n+1}z,\lambda^{n+1}t)-h_j(\lambda^nz,\lambda^nt)\right\} =$$
$${1\over {2i\pi }}\int_{\Bbb C}{{\left\{ (\psi_0(\xi )-1)+ (1-\psi_0(\xi ))\right\} f(\xi ,t)}\over {\xi -z}}d\xi \wedge d\overline \xi =0.$$
En r\'ealit\'e, la fonction ${\cal F}$-holomorphe $K$ qu'on cherche est r\'eduite \`a la s\'erie double :
$$K(z,t)=-\sum_{n=0}^\infty \sum_{j=1}^\infty \left\{ v_j(\lambda^{n+1}z,\lambda^{n+1}t)-v_j(\lambda^nz,\lambda^nt)\right\} .$$
Sa convergence uniforme sur tout compact de $\widetilde M$   r\'esulte de sa convergence en $(0,0)$
(puisque tous ses termes sont nuls)  et du fait que la s\'erie des d\'eriv\'ees par rapport
\`a la variable $z$ est \'equivalente \`a des s\'eries g\'eom\'etriques de raison   $\lambda $ (qui est dans $]0,1[$).
(Rappelons que sur l'espace vectoriel ${\cal H}_{\cal F}(\widetilde M)={\cal H}_{\cal F}({\Bbb C}\times {\Bbb R})$ la convergence
uniforme sur tout compact est \'equivalente \`a la convergence pour la $C^{0,\infty }$-topologie.)
\smallskip
\underbar{Sixi\`eme \'etape}
\smallskip
\noindent Revenons \`a la condition  $H(0,0)=0$ n\'ecessaire \`a l'existence de la fonction ${\cal F}$-holomorphe
$K$ qu'on cherche. Soit ${\cal I}:A_{\cal F}^{01}(M)\longrightarrow {\Bbb C}$ la forme lin\'eaire continue d\'efinie par :
$${\cal I}(fd\overline z)= {1\over {2i\pi }}\int_{\Bbb C}{{\psi_1(z)f(z,0)}\over {z}}dz\wedge d\overline z .$$
Alors il n'est pas difficile de voir, \`a partir de tous les calculs que nous avons men\'es pr\'ec\'edemment,
que $H(0,0)=0$ si, et seulement si, la $(0,1)$-forme feuillet\'ee $f(z,t)d\overline z$ est dans le noyau de
${\cal I}$.
La dimension de l'espace vectoriel $H_{\cal F}^{01}(M)$ est donc au plus $1$. Pour montrer   qu'elle est en fait \'egale \`a $1$,
il suffit de v\'erifier que la forme lin\'eaire   ${\cal I}$ est non nulle. Nous allons voir que son \'evaluation sur la  $(0,1)$-forme :
$$\omega_0={{zd\overline z}\over {z\overline z+t^2}}$$
est diff\'erente de $0$. \`A cet effet, rappelons d'abord que $\psi_1$ ne d\'epend que du module
de $z$ et que son support est contenu dans une couronne :
$$\Delta = \{ z\in {\Bbb C}:R_1 \leq \vert z\vert \leq R_2\} $$ (avec, bien s\^ur,
$0<R_1<R_2$). D'autre part, comme $\psi_1$ est \`a valeurs dans ${\Bbb R}_+$ et non identiquement nulle, son int\'egrale
sur l'intervalle $[R_1,R_2]$ (en tant que fonction uniquement de $r=\vert z\vert $) est un r\'eel strictement positif.
Posons  $z=re^{i\theta }$. On a :
$$dz\wedge d\overline z=d(re^{i\theta })\wedge d(re^{-i\theta })=\left( e^{i\theta }(dr+ird\theta )\right) \wedge
\left( e^{-i\theta }(dr-ird\theta )\right) =-2irdr\wedge d\theta .$$
D'o\`u :
$$\eqalign{{\cal I}(\omega_0)&= {1\over {2i\pi }}\int_{\Bbb C}{{\psi_1(z)z}\over {z\cdot \vert z\vert ^2}}dz\wedge d\overline z\cr
&=
{1\over {2i\pi }}\int_\Delta {{\psi_1(r)(-2i)}\over {r}}dr\wedge d\theta \cr
& =-{1\over \pi }\left( \int_{R_1}^{R_2}{{\psi_1(r)}\over r}dr\right) \cdot \int_0^{2\pi }d\theta \cr
&= -2\int_{R_1}^{R_2}{{\psi_1(r)}\over r}dr\cr
&<0.}$$
Par suite, la forme lin\'eaire continue ${\cal I}$ n'est pas nulle.
On a finalement :
$$H_{\cal F}^{01}({\Bbb S}^{n+1}\times {\Bbb S}^1)={\Bbb C}\omega_0.$$
Ceci   termine la d\'emonstration du th\'eor\`eme. \hfill $\square $
\smallskip
\noindent {\tfpp 4.4. Remarque}
\smallskip
Elle pourrait \^etre significative, et c'est la raison pour laquelle nous avons jug\'e de la faire. La cohomologie de Dolbeault
feuillet\'ee $H_{\cal F}^{0\ast }({\Bbb S}^{n+1}\times {\Bbb S}^1)$ du feuilletage ${\cal F}$ est la ``m\^eme" que celle de la feuille compacte
(courbe elliptique $C_\lambda $) induite par la feuille correspondant \`a $t=0$ dans le rev\^etement $\widetilde{M}={\Bbb C}\times {\Bbb R}\setminus\{ (0,0)\} $.
Ceci est s\^urement d\^u au fait que cette feuille compacte a une holonomie contractante.

\vskip0.5cm

\centerline{\tfp R\'ef\'erences}
\bigskip

\item{[Cam]} {\sc Camacho, C. \& Neto, A.} {\it Geometric Theory
of Foliations.} Birkh\"auser.

\smallskip

\item{[BC]} {\sc  Ben Charrada, R.}  {\it Cohomology of some complex laminations.}  Results in Math. 57, (2010) 33-41.

\smallskip

\item{[Ek1]} {\sc El Kacimi Alaoui, A.}  {\it The $\overline{\partial }$ along the leaves
and Guichard's Theorem for a simple complex foliation.}
Math. Annalen 347, (2010), 885-897.
\smallskip

\item{[Ek2]} {\sc El Kacimi Alaoui, A.}  {\it On leafwise meromorphic functions with prescribed poles.}
Bulletin of the Brazilian Mathematical Society, New Series, 48(2), (2017) 261-282.
\smallskip

\item{[ES]} {\sc El Kacimi Alaoui, A. \& Slim\`ene, J.}   {\it Cohomologie de Dolbeault le long des feuilles de
certains feuilletages complexes.}
Annales de l'Institut Fourier de Grenoble, Tome 60 n$^o$2, (2010), 727-757.
\smallskip

\item{[GT]} {\sc  Gigante, G. \& Tomassini, G.}  {\it Foliations with
complex leaves.}  Diff. Geo. and its Applications 5, (1995) 33-49.

\smallskip
\item{[H\"o]} {\sc H\"ormander, L.} {\it An Introduction to Complex
Analysis in Several Variables.}    D. Van Nostrand Compagny, Inc.,
(1966).

\smallskip
\item{[S$\ell $]} {\sc Slim\`ene, J.} {\it Deux exemples de calcul
explicite de cohomologie de Dolbeault feuillet\'ee}.
Proyecciones Vol. 27, N$^0$ 1, pp. 63-80, May 2008.

\bigskip
\bigskip
\bigskip

\noindent Rochdi {\sc Ben Charrada}

\noindent D\'epartement de Math\'ematiques

\noindent Facult\'e des Sciences de Sfax

\noindent 3018 Sfax -- Tunisie

\smallskip

\noindent rochdi$_-$charrada@yahoo.fr

\bigskip

\noindent Aziz {\sc El Kacimi Alaoui}

\noindent {\pg Université Polytechnique Hauts-de-France

\noindent EA 4015 - LAMAV

\noindent FR CNRS 2956

\noindent F-59313 Valenciennes Cedex 9, France}

\smallskip
\noindent aziz.elkacimi@uphf.fr

\end